\documentclass[reqno]{amsart}

\usepackage{amsmath,amsthm,amssymb,bm}
\usepackage{hyperref}
\usepackage{a4wide}
\usepackage{cleveref}
\usepackage{graphicx,color}
\usepackage{tikz}
\numberwithin{equation}{section}

\newtheorem{introthm}{Theorem}

\newtheorem{thm}{Theorem}[section]
\newtheorem{lem}[thm]{Lemma}
\newtheorem{prop}[thm]{Proposition}
\newtheorem{cor}[thm]{Corollary}
\theoremstyle{definition}
\newtheorem{exam}[thm]{Example}
\newtheorem{defn}[thm]{Definition}
\newtheorem{conj}[thm]{Conjecture}

\newtheorem{problem}[thm]{Problem}
\newtheorem{remark}[thm]{Remark}

\crefname{lem}{Lemma}{Lemmas}
\crefname{thm}{Theorem}{Theorems}
\crefname{prop}{Proposition}{Propositions}
\crefname{question}{Question}{Questions}
\crefname{defn}{Definition}{Definitions}
\crefname{conj}{Conjecture}{Conjectures}
\crefname{figure}{Figure}{Figures}
\crefname{cor}{Corollary}{Corollaries} 
\crefformat{equation}{(#2#1#3)}
\Crefformat{equation}{Equation #2(#1)#3}

\AddToHook{env/lem/begin}{\crefalias{thm}{lem}}
\AddToHook{env/prop/begin}{\crefalias{thm}{prop}}
\AddToHook{env/cor/begin}{\crefalias{thm}{cor}}
\AddToHook{env/exam/begin}{\crefalias{thm}{exam}}
\AddToHook{env/defn/begin}{\crefalias{thm}{defn}}
\AddToHook{env/conj/begin}{\crefalias{thm}{conj}}
\AddToHook{env/question/begin}{\crefalias{thm}{question}}
\AddToHook{env/problem/begin}{\crefalias{thm}{problem}}
\AddToHook{env/remark/begin}{\crefalias{thm}{remark}}

\newcommand\x{\mathbf{x}}
\newcommand\wt{\operatorname{wt}}
\newcommand\cro{\operatorname{cr}}
\newcommand\sgn{\operatorname{sgn}}
\newcommand\LL{\mathcal{L}}
\newcommand\Par{\operatorname{Par}}

\title[Positivity of Askey--Wilson moments]{Positivity and multivariate extensions of two-parameter Askey--Wilson moments}
\author{Younggwang Cho}
\thanks{}
\address{Department of Mathematics,
  Sungkyunkwan University (SKKU), Suwon, Gyeonggi-do 16419, South Korea}
\email{brglory@g.skku.edu}

\begin{document}

\begin{abstract}
  Kim and Stanton showed that certain specialized Askey--Wilson moments
  are polynomials with integer coefficients after normalization,
  and conjectured that these coefficients are nonnegative.
  We refine this into the conjecture for
  the normalized Askey--Wilson mixed moments and prove the refinement
  by establishing a manifestly positive recurrence.
  We also construct a combinatorial model that expresses the normalized
  Askey--Wilson mixed moments as weighted sums of decorated partial matchings.
  This model specializes to the partial matching model for the mixed moments of
  continuous \( q \)-Hermite polynomials.
  Finally, we consider the corresponding multivariate moments and prove that
  their minimal numerators have nonnegative coefficients up to an explicit sign.
\end{abstract}

\maketitle


\section{Introduction}
Orthogonal polynomials have been studied throughout
broad areas of mathematics, including
probability~\cite{KarlinMcGregor1957,Johansson2000,BrycWesolowski2010},
analysis~\cite{Deift1999,Lubinsky2009},
algebra~\cite{Macdonald2000,Cherednik1995,Koornwinder1992}
and combinatorics~\cite{Flajolet1980,Viennot1985,ISV}.
In particular, three-term recurrences have a crucial role in interpreting
quantities related to orthogonal polynomials in terms of
weighted lattice path models~\cite{Flajolet1980,Viennot1985}.

For an orthogonality functional \( \LL \) with \( \LL(1) = 1 \),
the moments \( \mu_n = \LL(x^n) \) determine its values on
all polynomials and hence determine the corresponding monic orthogonal polynomials~\cite[Chapter~I]{Chihara1978}.
Writing this sequence as \( \{p_n(x)\}_{n \ge 0} \),
we define the mixed moments by the change of basis from
monomials to given orthogonal polynomials:
\[
  x^n = \sum_{k \ge 0} \mu_{n, k} p_k(x).
\]
In particular, we have \( \mu_{n, 0} = \mu_n \).
Combinatorial formulas for moments or mixed moments provide
richer properties, such as coefficientwise nonnegativity,
that are not apparent from the integral representations or explicit algebraic formulas.

The Askey--Wilson polynomials \( p^{AW}_n(x; a, b, c, d; q) \)
lie at the top of the Askey scheme hierarchy and specialize to
other important orthogonal polynomial families in the scheme~\cite{AskeyWilson1985,Koekoek1996}.
In particular, their moments \( \mu_n(a, b, c, d; q) \)
have been studied from a combinatorial perspective~\cite{ISV,KS14}.
These moments have combinatorial interpretations involving
staircase tableaux~\cite{Corteel2011,Corteel2012} and
lecture hall graphs~\cite{Corteel2023}.
They are also connected with the asymmetric simple exclusion process
(ASEP)~\cite{Uchiyama2004,Corteel2011}.
The staircase tableau formulas are written in terms of the ASEP parameters
\( (\alpha, \beta, \gamma, \delta) \), which are related to
\( (a, b, c, d) \) by reparametrization~\cite{Corteel2011,Corteel2012}.

Although explicit formulas for the Askey--Wilson moments are known,
they do not directly show the nonnegativity of moments.
Even when an integral is nonnegative for a range of parameter values,
its expansion can have negative coefficients.
Conversely, a polynomial with nonnegative coefficients is nonnegative
whenever the variables are nonnegative.
In fact, the coefficientwise nonnegativity does depend on the choice of
parameters.
For example, the combinatorial formula involving
staircase tableaux~\cite{Corteel2011,Corteel2012}
shows the coefficientwise nonnegativity of the rescaled
Askey--Wilson moments in ASEP parameters
\( (\alpha, \beta, \gamma, \delta) \),
while this does not imply the nonnegativity in \( (a, b, c, d) \).

We study the two-parameter Askey--Wilson polynomials
obtained by the specialization \( b = -a \) and \( d = -c \).
Ismail and Simeonov obtained connection relations and a
characterization of this two-parameter Askey--Wilson family~\cite[Section~5.1]{IsmailSimeonov2012}.
Bryc and Weso{\l}owski used a purely imaginary specialization
of this family to construct purely
quadratic harnesses~\cite[Section~5.4]{BrycWesolowski2010}.

Within this family, \( c^2 = a^2q \) gives, up to normalization,
the continuous \( q \)-ultraspherical polynomials
(or Rogers polynomials)~\cite[Section~3.10.1]{Koekoek1996}.
The mixed moments of Rogers polynomials were obtained by
Ismail and Rahman~\cite[Theorem~4.2]{IsmailRahman2011},
and recent work of Zeng studies Rogers-moment deformations related to
Narayana polynomials and Hankel positivity~\cite{Zeng2026Rogers}.
The further specialization \( a = c = 0 \) gives the continuous
\( q \)-Hermite family.

For \( b = -a \), \( d = -c \), the odd moments vanish,
and the even moments are rational functions of
\( A = a^2, C = c^2 \) and \( q \).
We write 
\[
  (z; q)_r = \prod_{j = 0}^{r-1}(1 - zq^j)
\]
for the \( q \)-Pochhammer symbol, with the convention \( (z; q)_0 = 1 \).
Following Kim and Stanton~\cite[Proposition~2.9 and Conjecture~1]{KS14}, we consider
the following normalization of the even Askey--Wilson moments:
\[
  \tau_{2n}(A, C; q) = \tau_{2n}(a^2, c^2; q)
  = \left(\frac{4}{1-q}\right)^{n}(a^2c^2q ; q^2)_n \mu_{2 n}(a,-a, c,-c ; q).
\]
They showed that the normalized Askey--Wilson moment
\( \tau_{2n}(A, C; q) \) is a polynomial in \( A, C, q \)
and conjectured that \( \tau_{2n}(A, C; q) \) has
nonnegative integer coefficients.

Our first result extends this conjecture to mixed moments.
Let \( \{p^\nu_n(x; A, C; q)\}_{n \ge 0} \)
be the orthogonal polynomial sequence satisfying the recurrence
\[
  p^\nu_{n+1} = xp^\nu_n - \lambda^\nu_n p^\nu_{n-1}
\]
where
\[
  \lambda^\nu_n = \lambda^\nu_n(A, C; q)
  = \frac{(1+Aq^{n-1})(1+Cq^{n-1})(1-ACq^{n-2})}
  {(1-ACq^{2n-3})(1-ACq^{2n-1})}[n]_q.
\]
Here and throughout the paper, \( [n]_q = 1 + q + \cdots + q^{n-1} \).
Then, the Askey--Wilson mixed moments \( \mu_{n, k}(a, b, c, d; q) \) and
the rescaled mixed moments \( \nu_{n, k}(A, C; q) \) are related by
\[
  \nu_{n, k}(a^2, c^2; q) = \left(\frac{2}{\sqrt{1-q}}\right)^{n-k}
  \mu_{n, k}(a, -a, c, -c; q).
\]
For \( 0 \le k \le n \) and \( n - k \) even,
define the normalized mixed moments
\( \tau_{n, k}(A, C; q) \) to be
\[
  \tau_{n, k}(A, C; q) = 
  (ACq^{2k+1}; q^2)_{\frac{n-k}{2}} \nu_{n, k}(A, C; q).
\]
In particular we have \( \tau_{2n , 0}(A, C; q) = \tau_{2n}(A, C; q) \).

We prove the coefficientwise nonnegativity of  \( \tau_{n, k}(A, C; q) \)
through a positive recurrence.
When \( k = 0 \), this result proves the Kim--Stanton conjecture.
\begin{introthm}[\Cref{thm:main}]
  Let \( 0 \le k \le n \) with \( n - k \) even.
  Then \( \tau_{n,k}(A,C;q) \) satisfies the following recurrence:
  \[
    \tau_{n, k}(A, C; q) = \tau_{n-1, k-1}(A, Cq^2; q)
    + [k+1]_q(1+Aq^k)(1+C)\tau_{n-1, k+1}(A, C; q) \quad \text{for }n \ge 1,
  \]
  with boundary conditions
  \[
    \tau_{0, 0}(A, C; q) = 1, \qquad \tau_{n, k} = 0
    \text{ unless } 0 \le k \le n \text{ and } n-k \text{ is even}.
  \]
\end{introthm}

Our second result is a combinatorial model
for the normalized Askey--Wilson mixed moments \( \tau_{n, k}(A, C; q) \).
By \cite[Proposition~2.9]{KS14}, the sum of coefficients is
\[
  \tau_{2n}(1, 1; 1) = 4^n(2n-1)(2n-3) \cdots 1.
\]
Since the number of perfect matchings on \( [2n] \) is \( (2n-1)(2n-3) \cdots 1 \),
this identity suggests a model for \( \tau_{2n}(A, C; q) \)
based on perfect matchings.

Furthermore, the fact that the sum of coefficients of \( \tau_{n, k}(A, C; q) \) is
\[
  \tau_{n, k}(1, 1; 1) =
  \begin{cases}
    2^{n-k}\binom{n}{k}(n-k-1)(n-k-3) \cdots 1 & n-k \text{ is even}, \\
    0 & \text{otherwise},
  \end{cases}
\]
suggests a model for \( \tau_{n, k}(A, C; q) \) based on partial matchings.
Moreover, the factor \( 2^{n-k} \) suggests that each nonfixed vertex may be
independently decorated or left undecorated.
Motivated by this, we construct a weighted model of decorated matchings of \( [n] \)
with \( k \) fixed points.
We denote this set by \( M^*(n, k) \),
and assign a weight to each matching.
For details, see \Cref{sec:com-model}.
\begin{introthm}[\Cref{thm:comb}]
  We have
  \[
    \tau_{n, k}(A, C; q) = \sum_{\sigma \in M^*(n, k)} \wt(\sigma).
  \]
\end{introthm}
Under the specialization \( A = C = 0 \),
this model reduces to the partial matching model for the mixed moments
of continuous \( q \)-Hermite polynomials~\cite{CiglerZeng2011}.

We next turn to multivariate moments.
The multivariate moments of Askey--Wilson polynomials
are motivated by the relation between Koornwinder moments under \( t = q \)
and the two-species ASEP~\cite{Corteel2017,Corteel2018}.
In that setting, Rains conjectured the coefficientwise nonnegativity
of normalized multivariate moments in the ASEP parameters
\( (\alpha, \beta, \gamma, \delta) \), \( q \) and \( \xi \)~\cite[Conjecture~4.4]{Corteel2018}.
Special cases of this conjecture have been established in \cite{Cho2026a,Cho2026c}.
We consider the corresponding question for
the two-parameter Askey--Wilson polynomials.

Fix \( n \ge 0 \) and let \( \Par_n \) be the set of
weakly decreasing sequences of nonnegative integers of length \( n \).
Write \( \emptyset = (0 ,\dots, 0) \), and let \( \x_n = (x_1 ,\dots, x_n) \)
be a tuple of \( n \) variables.
The multivariate lift of \( p_n^\nu(x; A, C; q) \) is
\[
  p^\nu_\lambda(\x_n; A, C; q) =
  \frac{\det(p^\nu_{\lambda_i + n - i}(x_j))_{i,j = 1}^n}
  {\det(x_j^{n-i})_{i, j = 1}^n}.
\]
Define the Schur polynomial
\[
  s_\lambda(\x_n) := \frac{\det(x_j^{\lambda_i + n - i})_{i, j = 1}^n}
  {\det(x_j^{n-i})_{i, j = 1}^n}.
\]
Then the multivariate mixed moments \( V_{\lambda, \rho}(A, C; q) \) are defined by
\[
  s_\lambda(\x_n) = \sum_{\rho \in \Par_n} V_{\lambda, \rho}(A, C; q)p^\nu_\rho(\x_n; A, C; q).
\]
In particular, \( V_\lambda(A, C; q) := V_{\lambda, \emptyset}(A, C; q) \)
is called the multivariate moment.
The multivariate moment has the determinantal formula~\cite[Proposition~2.10]{LHT}:
\[
  V_\lambda(A, C; q) = \det(\nu_{\lambda_i+n-i, n-j})_{i,j=1}^n.
\]
Let \( T_\lambda(A, C; q) \) be the minimal numerator of \( V_\lambda(A, C; q) \).

Our third main result is the nonnegativity of \( T_\lambda(A, C; q) \)
up to an explicit sign.
The nonnegativity of \( \tau_{n, k}(A, C; q) \) does not immediately imply
the nonnegativity of \( T_{\lambda}(A, C; q) \).
Thus, we construct a lattice path model for another determinantal quantity
\( T'_\lambda(A, C; q) \) and use the Lindstr\"{o}m--Gessel--Viennot lemma to prove that \( T'_\lambda(A, C; q) \)
has nonnegative coefficients up to an explicit sign.
We then show that \( T_{\lambda}(A, C; q) = T'_\lambda(A, C; q) \),
which proves the nonnegativity of \( T_{\lambda}(A, C; q) \) up to sign.
\begin{introthm}[\Cref{thm:multi}]
  Let \( \lambda = (\lambda_1 ,\dots, \lambda_n) \in \Par_n \). Then
  \[
    (-1)^{\lambda_2 + \lambda_4 + \cdots } T_{\lambda}(A, C; q)
  \]
  is a polynomial in \( A, C, q \) with nonnegative integer coefficients.
\end{introthm}

The paper is organized as follows.
In \Cref{sec:prelim}, we review the basic properties of orthogonal polynomials
and the Kim--Stanton conjecture about the normalized moments of Askey--Wilson polynomials.
In \Cref{sec:proof-of-KSconj}, we formulate and prove a refinement of
the Kim--Stanton conjecture.
In \Cref{sec:com-model}, we construct a combinatorial model
based on the theorem proved in \Cref{sec:proof-of-KSconj}.
In \Cref{sec:multi}, we extend the result to multivariate moments
and use a lattice path model to prove that their minimal numerators
have nonnegative coefficients up to sign.

\section{Preliminaries} \label{sec:prelim}
\subsection{Basic definitions}
A \emph{partition} is a weakly decreasing sequence
\( \lambda = (\lambda_1 \ge \cdots \ge \lambda_n) \) of
nonnegative integers.
The sum \( \lambda_1 + \cdots + \lambda_n \) is denoted by \( |\lambda| \).
The collection of all partitions having \( n \) nonnegative parts is denoted by
\( \Par_n \), that is,
\[
  \Par_n = \{ (\lambda_1 ,\dots, \lambda_n) \in \mathbb{Z}_{\ge0}^n :
  \lambda_1 \ge \cdots \ge \lambda_n \}.
\]
For \( \lambda, \rho \in \Par_n \),
we write \( \rho \subseteq \lambda \)
if \( \rho_i \le \lambda_i \) for all \( 1 \le i \le n \).
We denote \( (0 ,\dots, 0) \in \Par_n \) by \( \emptyset \).

We use the \( q \)-Pochhammer symbols:
\[
  (a;q)_n := \prod_{k=0}^{n-1}(1-aq^k),
  \qquad (a_1, \dots, a_m; q)_{n} := (a_1; q)_{n} \cdots (a_m; q)_{n},
\]
where \( n \) may be \( \infty \), in which case the product is over all \( k\ge0 \).

The \( q \)-integer \( [n]_q \) is defined by
\[
  [n]_q := 1 + \cdots + q^{n-1}.
\]

The set \( \{ 1, 2 ,\dots, n \} \) is simply denoted by \( [n] \).

The \emph{symmetric group} \( S_n \) is the group of bijections on \( [n] \),
equipped with the composition of functions as the multiplication.
An \emph{involution} in the symmetric group is a bijection
whose square is the identity.

A \emph{Dyck path} is a lattice path consisting of
up steps \( U = (1, 1) \) and down steps \( D = (1, -1) \)
which does not go below the \( x \)-axis.
Two paths are \emph{intersecting} if they share a vertex.
Each shared vertex is called an \emph{intersection}.

\subsection{Orthogonal polynomials}
In this subsection, we define orthogonal polynomials and
recall their basic properties.
Throughout this paper, all orthogonal polynomials are assumed to be monic.

Let \( \LL \) be a linear functional defined on
the space of polynomials in \( x \) with \( \LL(1) = 1 \).
Let \( \{p_n(x)\}_{n \ge 0} \) be a sequence of monic polynomials
with \( \deg(p_n(x)) = n \).
It is called an \emph{orthogonal polynomial sequence}
with respect to \( \LL \) if it satisfies
\[
  \LL(p_n(x)p_m(x)) = K_n\delta_{n, m}, \qquad K_n \neq 0,
\]
where \( \delta_{n, m} = 1 \) if \( n = m \), and \( 0 \) otherwise.
By Favard's theorem, a sequence of monic orthogonal polynomials
\( \{p_n(x)\}_{n \ge 0} \) satisfies a three-term recurrence
of the form
\[
  p_{n+1}(x) = (x - b_n)p_n(x) - \lambda_n p_{n-1}(x) \quad \text{for } n \ge 0,
  \quad p_{-1}(x) = 0, \quad p_0(x) = 1,
\]
with \( \lambda_n \neq 0 \) for \( n \ge 1 \).

For a monic orthogonal polynomial sequence \( \{p_n(x)\}_{n \ge 0} \),
the \emph{mixed moments} \( \mu_{n, k} \) are defined by
\[
  x^n = \sum_{k = 0}^n \mu_{n, k} p_k(x).
\]
For \( k = 0 \), these are called \emph{moments} and denoted by \( \mu_n = \mu_{n, 0} \).

We now construct the multivariate lift of a monic orthogonal polynomial sequence.
Let \( \x_n = (x_1 ,\dots, x_n) \) be a tuple of variables.
For \( \lambda \in \Par_n \), the \emph{Schur polynomial} \( s_\lambda(\x_n) \)
is defined by
\[
  s_\lambda(\x_n) := \frac{\det(x_j^{\lambda_i + n - i})_{i, j = 1}^n}
  {\Delta(\x_n)},
  \qquad \Delta(\x_n) :=  \det(x_j^{n-i})_{i, j = 1}^n.
\]
Similarly, for a monic orthogonal polynomial sequence \( \{p_n(x)\}_{n \ge0 } \),
the corresponding \emph{multivariate orthogonal polynomials}
indexed by partitions \( \lambda \in \Par_n \) are defined by
the following bialternant formula:
\[
  p_\lambda(\x_n) = \frac{\det(p_{\lambda_i + n - i}(x_j))_{i,j = 1}^n}
  {\Delta(\x_n)}.
\]
For \( \lambda, \rho \in \Par_n \),
the \emph{multivariate mixed moments} \( M_{\lambda, \rho} \) are defined by
\[
  s_\lambda(\x_n) = \sum_{\rho \in \Par_n} M_{\lambda, \rho} p_\rho(\x_n).
\]
For \( \rho = \emptyset \), these are called \emph{multivariate moments}
and denoted by \( M_\lambda = M_{\lambda, \emptyset} \).
The multivariate mixed moments \( M_{\lambda, \rho} \) and
the univariate mixed moments \( \mu_{n, k} \) are related by the following
determinantal formula~\cite[Proposition~2.10]{LHT}:
\[
  M_{\lambda, \rho} = \det(\mu_{\lambda_i+n-i, \rho_j+n-j})_{i,j=1}^n.
\]

\subsection{Askey--Wilson polynomials and the Kim--Stanton conjecture}
In this subsection, we recall the definition of Askey--Wilson polynomials
in terms of their three-term recurrence and consider the specialization \( b = -a, d = -c \).
We then review the normalization and positivity conjecture of Kim and Stanton.

The \emph{monic Askey--Wilson polynomials}
\( p^{AW}_{n}(x) = p^{AW}_n(x; a, b, c, d; q) \)
satisfy the three-term recurrence
\[
  p^{AW}_{n+1}(x) = \left(x - \frac{a + a^{-1} - A_n - C_n}{2}\right)
  p^{AW}_n(x) - \frac{1}{4}A_{n-1}C_n p^{AW}_{n-1}(x),
\]
where
\begin{align*}
  A_n &= \frac{(1-abq^n)(1-acq^n)(1-adq^n)(1-abcdq^{n-1})}{a(1-abcdq^{2n-1})(1-abcdq^{2n})}, \\
  C_n &= \frac{a(1-q^n)(1-bcq^{n-1})(1-bdq^{n-1})(1-cdq^{n-1})}{(1-abcdq^{2n-2})(1-abcdq^{2n-1})}.
\end{align*}
Under the specialization \( b=-a \) and \( d=-c \),
the three-term recurrence reduces to
\[
  p^{AW}_{n+1}(x; a, -a, c, -c; q) = x p^{AW}_{n}(x; a, -a, c, -c; q)
  - \lambda_n p^{AW}_{n-1}(x; a, -a, c, -c; q),
\]
where
\[
  \lambda_n =
  \frac{(1-q^n)(1+a^2q^{n-1})(1+c^2q^{n-1})(1-a^2c^2q^{n-2})}
  {4(1-a^2c^2q^{2n-3})(1-a^2c^2q^{2n-1})}.
\]
Since the specialized polynomials and moments depend only on \( a^2 \) and \( c^2 \),
we write \( A = a^2 \) and \( C = c^2 \) throughout the paper.

Following Kim and Stanton~\cite{KS14}, we define the
\emph{normalized Askey--Wilson moment} \( \tau_{2n}(A, C; q) \)
of the two-parameter Askey--Wilson polynomials by
\[
  \tau_{2n}(A, C; q) = \tau_{2 n}(a^2, c^2; q)=
  4^n(a^2c^2q ; q^2)_n \mu_{2 n}(a,-a, c,-c ; q) /(1-q)^n.
\]
\begin{prop}\cite[Proposition~2.9]{KS14} \label{prop:2.9}
  The normalized Askey--Wilson moments \( \tau_{2n}(A, C; q) \)
  are polynomials in \( A, C, q \) with integer coefficients.
  The sums of coefficients are given by
  \[
    \tau_{2n}(1, 1; 1) = 2^{2 n}(2 n-1)(2 n-3) \cdots 1.
  \]
\end{prop}
They further conjectured that the normalized Askey--Wilson moments
have nonnegative coefficients.
\begin{conj} \label{conj:main} \cite[Conjecture~1]{KS14}
  The coefficients of \( \tau_{2n}(A, C; q) \) are nonnegative integers.
\end{conj}

\section{A refinement and its proof}
\label{sec:proof-of-KSconj}
In this section, we refine \Cref{conj:main} by considering the nonnegativity of
the normalized mixed moments.
We construct a recurrence which preserves nonnegativity.
Since the boundary values are nonnegative,
the nonnegativity of the normalized mixed moments follows by induction.

In the proof of \Cref{prop:2.9}, Kim and Stanton introduce
a rescaled version of the two-parameter Askey--Wilson polynomials
\( p_n^\nu=p_n^\nu(x; A, C; q) \)
satisfying the three-term recurrence
\begin{equation}\label{eq:recurrence-nu}
  p^\nu_{n+1} = xp^\nu_n - \lambda^\nu_n p^\nu_{n-1},
\end{equation}
where
\[
  \lambda^\nu_n = \lambda^\nu_n(A, C; q)
  = \frac{(1+Aq^{n-1})(1+Cq^{n-1})(1-ACq^{n-2})}
  {(1-ACq^{2n-3})(1-ACq^{2n-1})}[n]_q.
\]

Let \( \nu_{2n}(A, C; q) \) denote the moments of the polynomial sequence
\( \{p_n^\nu(x; A, C; q)\}_{n \ge 0} \).
Comparing the two three-term recurrences, we obtain
\[
  \nu_{n, k}(A, C; q)
  = \left( \frac{2}{\sqrt{1-q}} \right)^{n-k} \mu_{n, k}(a, -a, c, -c; q).
\]
In particular,
\[
  \tau_{2n}(A, C; q) = (ACq; q^2)_n \nu_{2n}(A, C; q).
\]

Now we fix our convention for numerators and denominators.
For a rational function \( f(A, C, q) \in\mathbb{Q}(A,C,q) \),
write \( f(A, C, q)=P(A, C, q)/Q(A, C, q) \),
where \( P(A,C,q), Q(A,C,q)\in\mathbb{Q}[A,C,q] \) are coprime.
Furthermore, if \( Q(0, 0, 0) \neq 0 \),
it is always possible to rescale \( P(A, C, q) \) and \( Q(A, C, q) \)
so that \( Q(0,0,0)=1 \).
We call such \( P(A, C, q) \) and \( Q(A, C, q) \) the \emph{minimal numerator} and
the \emph{minimal denominator} of \( f(A, C, q) \), respectively.
This normalization determines \( P(A, C, q) \) and \( Q(A, C, q) \) uniquely.

We now define the normalized Askey--Wilson mixed moments.
\begin{defn}\label{def:normalized-mixed}
  The \emph{normalized mixed moment} \( \tau_{n, k}(A, C; q) \) is
  defined by
  \[
    \tau_{n, k}(A, C; q) =
    \begin{cases}
      (ACq^{2k+1}; q^2)_{\frac{n-k}{2}} \nu_{n, k}(A, C; q)
      & 0 \le k \le n, n - k \mbox{ is even,}\\
      0 & \mbox{otherwise.}
    \end{cases}
  \]
\end{defn}

For \( 0\le k\le n \) with \( n-k \) even, the normalizing factor
\( (ACq^{2k+1};q^2)_{(n-k)/2} \) is precisely the minimal denominator
of \( \nu_{n,k}(A,C;q) \).
Thus, \( \tau_{n,k}(A,C;q) \) is its minimal numerator.
We prove this in \Cref{cor:min-deno} as a consequence of \Cref{thm:main}.

We now state the refined result.
When \( k = 0 \), this result proves the Kim--Stanton conjecture.
\begin{thm} \label{thm:main}
  Let \( 0 \le k \le n \) with \( n - k \) even.
  Then \( \tau_{n,k}(A,C;q) \) satisfies the following recurrence:
  \begin{equation}\label{eq:main}
    \tau_{n, k}(A, C; q) = \tau_{n-1, k-1}(A, Cq^2; q)
    + [k+1]_q(1+Aq^k)(1+C)\tau_{n-1, k+1}(A, C; q) \quad \text{for }n \ge 1,
  \end{equation}
  with boundary conditions
  \[
    \tau_{0, 0}(A, C; q) = 1, \qquad \tau_{n, k} = 0
    \text{ unless } 0 \le k \le n \text{ and } n-k \text{ is even}.
  \]

  In particular, \( \tau_{n, k}(A, C; q) \) is a polynomial
  in \( A, C, q \) with nonnegative coefficients.
\end{thm}

Before proving the theorem, we introduce the notation
\[
  p_n = p^\nu_n(x; A, C; q), \quad p'_n = p^\nu_n(x; A, Cq^2; q),
  \quad \lambda_n = \lambda^\nu_n(A, C; q),
  \quad \lambda'_n = \lambda^\nu_n(A, Cq^2; q).
\]
We use the convention that
\[
  p_{n} = 0, \quad p'_n = 0, \quad \lambda_n = 0, \quad \lambda'_n = 0
\]
for \( n < 0 \), and
we set \( \nu_{n,k}=\nu'_{n,k}=0 \) unless \( 0 \le k \le n \).

Define 
\[
  \eta_n = \frac{C(1+Aq^n)(1+Aq^{n+1})(1-q)[n+1]_q [n+2]_q}{(1-ACq^{2n+1})(1-ACq^{2n+3})}.
\]
We set \( \eta_{n} = 0 \) if \( n < 0 \).
The following lemma relates
\( p_n \) and \( p'_n \).
\begin{lem} \label{lem:1}
  For \( n \ge 0 \), we have
  \[
    p_n = p'_n - \eta_{n-2}p'_{n-2}.
  \]
\end{lem}
\begin{proof}
  Before proving the lemma, we first establish the following two identities:
  For \( m \ge 2 \) we have
  \begin{align}
    \lambda'_{m-1} - \lambda_{m-1}
    &= \eta_{m-3} - \eta_{m-2}, \label{eq:first}\\
    \lambda_{m-1} \eta_{m-4} &= \eta_{m-3}\lambda'_{m-3}. \label{eq:second}
  \end{align}

  To verify \eqref{eq:first}, we substitute the definitions and obtain
  \begin{align*}
    &\lambda'_{m-1} - \lambda_{m-1} \\
    &\quad= \frac{(1+Aq^{m-2})[m-1]_q}{1-ACq^{2m-3}}\left(
      \frac{(1+Cq^{m})(1-ACq^{m-1})}{(1-ACq^{2m-1})}
      - \frac{(1+Cq^{m-2})(1-ACq^{m-3})}{(1-ACq^{2m-5})} \right) \\
    &\quad= \frac{C(1+Aq^{m-2})[m-1]_q}{1-ACq^{2m-3}}\left(
      \frac{(1+Aq^{m-3})(1-q^{m-2})}{(1-ACq^{2m-5})}
      - \frac{(1+Aq^{m-1})(1-q^{m})}{(1-ACq^{2m-1})} \right) \\
    &\quad=\eta_{m-3} - \eta_{m-2}.
  \end{align*}
  We similarly obtain
  \begin{align*}
    \lambda_{m-1} \eta_{m-4} &= \frac{C(1-q)(1+Aq^{m-4})(1+Aq^{m-3})(1+Aq^{m-2})}
      {(1-ACq^{2m-7})(1-ACq^{2m-5})^2(1-ACq^{2m-3})} \\
    &\qquad \times (1+Cq^{m-2})(1-ACq^{m-3})[m-3]_q [m-2]_q[m-1]_q \\
    &=\lambda'_{m-3}\eta_{m-3}.
  \end{align*}

  We now prove the lemma using induction on \( n \).
  The cases \( n = 0 \) and \( n = 1 \) are immediate.
  Suppose that \( n \ge 2 \).  
  Using the induction hypothesis and the three-term recurrences, we obtain
  \begin{align*}
    p_n
    &= xp_{n-1} - \lambda_{n-1}p_{n-2} \\
    &= xp'_{n-1} - \eta_{n-3}xp'_{n-3} - \lambda_{n-1}p'_{n-2}
      + \lambda_{n-1}\eta_{n-4}p'_{n-4} \\
    &= (p'_n + \lambda'_{n-1}p'_{n-2}) - \eta_{n-3}(p'_{n-2} + \lambda'_{n-3}p'_{n-4})
      - \lambda_{n-1}p'_{n-2} + \lambda_{n-1}\eta_{n-4}p'_{n-4} \\
    &= p'_n + (\lambda'_{n-1} - \eta_{n-3} -\lambda_{n-1})p'_{n-2}
      + (\lambda_{n-1}\eta_{n-4} - \eta_{n-3}\lambda'_{n-3})p'_{n-4} \\
    &= p'_n - \eta_{n-2}p'_{n-2}.
  \end{align*}
  The last equality follows from \eqref{eq:first} and \eqref{eq:second}.
\end{proof}

We now prove \Cref{thm:main}.

\begin{proof}[Proof of \Cref{thm:main}]
  Recall that the mixed moments \( \nu_{n, k} \)
  are the coefficients in the expansion of \( x^n \)
  in the basis \( \{p_k\}_{k \ge 0} \).
  That is,
  \[
    x^n =  \sum_{k=0}^n \nu_{n, k} p_k, \qquad x^n =  \sum_{k=0}^n \nu'_{n, k} p'_k.
  \]
  By \Cref{lem:1}, we obtain
  \[
    x^n = \sum_{k=0}^n \nu_{n, k} (p'_k - \eta_{k-2}p'_{k-2})
        = \sum_{k=0}^n (\nu_{n, k} - \nu_{n, k+2}\eta_k)p'_k.
  \]
  Thus we have
  \( \nu'_{n, k} = \nu_{n, k} - \eta_k\nu_{n, k+2} \),
  or equivalently
  \begin{equation} \label{eq:1}
    \nu_{n, k} = \nu'_{n, k} + \eta_k\nu_{n, k+2}.
  \end{equation}

  From the definition of mixed moments and \Cref{eq:recurrence-nu}, we have
  \[
    \nu_{n+1, k} = \nu_{n, k-1} + \lambda_{k+1} \nu_{n, k+1}.
  \]
  By \eqref{eq:1}, we have
  \[
    \nu_{n+1, k} = (\nu'_{n, k-1} + \eta_{k-1}\nu_{n, k+1}) + \lambda_{k+1} \nu_{n, k+1}
    = \nu'_{n, k-1} + (\eta_{k-1} + \lambda_{k+1}) \nu_{n, k+1}.
  \]

  The coefficient \( \eta_{k-1} + \lambda_{k+1} \) simplifies to
  \begin{align*}
    &\eta_{k-1} + \lambda_{k+1} \\
    &\quad= \frac{C(1+Aq^{k-1})(1+Aq^{k})(1-q)[k]_q [k+1]_q}{(1-ACq^{2k-1})(1-ACq^{2k+1})}
      + \frac{(1+Aq^{k})(1+Cq^{k})(1-ACq^{k-1})}{(1-ACq^{2k-1})(1-ACq^{2k+1})}[k+1]_q \\
    &\quad= \frac{(1+Aq^{k})[k+1]_q}{(1-ACq^{2k-1})(1-ACq^{2k+1})}
      (C(1+Aq^{k-1})(1-q^k) + (1+Cq^{k})(1-ACq^{k-1})) \\
    &\quad= \frac{(1+Aq^{k})(1+C)[k+1]_q}{(1-ACq^{2k+1})}.
  \end{align*}
  Therefore, we obtain
  \[
    \nu_{n+1, k} = \nu'_{n, k-1} +
    \frac{(1+Aq^{k})(1+C)[k+1]_q}{(1-ACq^{2k+1})} \nu_{n, k+1}.
  \]
  Equation \eqref{eq:main} follows by multiplying
  both sides by \( (ACq^{2k+1};q^2)_{(n-k+1)/2} \).
\end{proof}

\begin{cor}\label{cor:min-deno}
  Let \( 0\le k\le n \) with \( n-k \) even.
  The minimal denominator of \( \nu_{n,k}(A,C;q) \) is
  \[
    (ACq^{2k+1};q^2)_{(n-k)/2},
  \]
  and its minimal numerator is \( \tau_{n,k}(A,C;q) \).
\end{cor}
\begin{proof}
  By \Cref{thm:main} and the initial condition \( \tau_{0,0}=1 \),
  induction shows that \( \tau_{n,k} \) is a nonzero polynomial
  with nonnegative coefficients throughout the stated range.
  Let \( D_{n, k} = (ACq^{2k+1};q^2)_{(n-k)/2} \).
  Since \( \nu_{n,k}=\tau_{n,k}/D_{n,k} \) and \( D_{n,k} \) is a product
  of irreducible factors of the form \( 1-ACq^i \), we have
  \begin{equation}\label{eq:cor-3.4}
    \tau_{n,k}(A,C;q)=P_{n,k}(A,C;q)\prod_{i\in I}(1-ACq^i)
  \end{equation}
  for some set \( I \), where \( P_{n,k} \) is the minimal numerator of \( \nu_{n,k} \).

  Set \( A=C=q=1 \).
  The left-hand side of \eqref{eq:cor-3.4} is positive,
  whereas the right-hand side of \eqref{eq:cor-3.4} vanishes if \( I \) is nonempty.
  Hence \( I \) is empty, so \( \tau_{n,k}=P_{n,k} \) and
  \( D_{n,k} \) is the minimal denominator of \( \nu_{n,k} \).
\end{proof}

Using \eqref{eq:main}, we can construct a natural Dyck path model.
For any Dyck path segment \( p \) ending at height \( k \), we define
\( \wt(p) \) recursively.
\begin{itemize}
\item If \( p \) is empty, \( \wt(p) = 1 \).
\item If the last step is an up step, \( p = p' U \),
  we set \( \wt(p) = \wt(p')|_{C = Cq^2} \).
\item If the last step is a down step, \( p = p'D \),
  we set \( \wt(p) = [k+1]_q(1+Aq^k)(1+C)\wt(p') \).
\end{itemize}

By \eqref{eq:main},
\begin{equation}\label{eq:path}
  \tau_{n, k}(A, C; q) = \sum_{p} \wt(p),
\end{equation}
where the sum ranges over the Dyck paths from \( (0, 0) \) to \( (n, k) \).

\begin{exam}
  The two Dyck paths from \( (0,0) \) to \( (4,0) \) have weights
  \[
    \begin{aligned}
      \wt(UDUD) &= (1+A)^2(1+C)(1+Cq^2), \\
      \wt(UUDD) &= [2]_q(1+A)(1+Aq)(1+C)^2.
    \end{aligned}
  \]
  For \( UDUD \), the third step \( U \) makes the substitution
  \( C\mapsto Cq^2 \) in the entire weight \( (1+A)(1+C) \)
  of the initial segment \( UD \).
  By \eqref{eq:path}, the sum of these two weights is
  \( \tau_{4,0}(A,C;q) \).
\end{exam}

Since the weights depend only on the sequence of heights,
translating each path by \( (-n, 0) \) preserves its weight.
Therefore, \( \tau_{n, k}(A, C; q) \) is also the weighted sum
of the Dyck paths from \( (-n, 0) \) to \( (0, k) \).

\section{Combinatorial model}
\label{sec:com-model}
In this section, we construct a weighted combinatorial model
for the normalized mixed moment \( \tau_{n, k}(A, C; q) \).
This model gives a combinatorial explanation for
the nonnegativity of \( \tau_{n, k}(A, C; q) \).

\begin{defn}
  A \emph{(partial) matching} of \( [n] \) with \( k \) fixed points
  is an involution \( \sigma \in S_n \)
  with exactly \( k \) fixed points.
  We denote the set of all such partial matchings by \( M(n, k) \).
  We set \( M(n, k) = \emptyset \) unless \( 0 \le k \le n \)
  and \( n - k \) is even.

  Each \( 2 \)-cycle \( (i, j) \) with \( i < j \)
  in the disjoint cycle decomposition of \( \sigma \)
  is called an \emph{edge} with endpoints \( i \) and \( j \).
  We call \( i \) the \emph{starting point}
  and \( j \) the \emph{ending point}.
  
  Two edges \( (i_1, j_1) \) and \( (i_2, j_2) \) form a \emph{crossing}
  if \( i_1 < i_2 < j_1 < j_2 \).
  A fixed point \( (i_1) \) and an edge \( (i_2, j_2) \) form a \emph{crossing}
  if \( i_2 < i_1 < j_2 \).
  The number of crossings of a matching \( \sigma \) is denoted by
  \( \cro(\sigma) \).
  
  A matching with no fixed points is called a \emph{perfect matching}.
\end{defn}

\begin{remark}
  The number of matchings of \( [n] \) with \( k \) fixed points is
  \( \binom{n}{k}(n-k-1)(n-k-3) \cdots 1 \), when \( n - k \) is even.
  In particular, the number of perfect matchings on \( [2n] \) is
  \( (2n-1)(2n-3) \cdots 1 \).
\end{remark}

To construct a combinatorial model for \( \tau_{n, k}(A, C; q) \),
we introduce decorated matchings.

\begin{defn}
  A \emph{decorated matching} of \( [n] \) with \( k \) fixed points
  is a partial matching in \( M(n, k) \)
  in which each endpoint of each edge is
  independently either decorated or undecorated.
  Fixed points are left undecorated.

  The set of decorated matchings of \( [n] \) with \( k \) fixed points
  is denoted by \( M^*(n, k) \).
\end{defn}

\begin{remark}
  For each edge in a matching, there are four possible decoration patterns.
  Therefore, for \( 0 \le k \le n \) with \( n-k \) even, the number of decorated matchings of
  \( [n] \) with \( k \) fixed points is
  \( 2^{n-k}\binom{n}{k}(n-k-1)(n-k-3) \cdots 1 \) and
  the number of decorated perfect matchings on \( [2n] \) is
  \( 2^{2n}(2n-1)(2n-3) \cdots 1 \).
\end{remark}

Given a decorated matching \( \sigma \),
we assign a weight to each vertex as follows.
\begin{itemize}
\item If a vertex \( i \) is undecorated, set \( \wt_\sigma(i) = 1 \).
\item If a vertex \( i \) is a decorated starting point,
  let \( t_i \) be the number of fixed points smaller than \( \sigma(i) \),
  plus the number of edges \( (a, b) \) satisfying \( a < \sigma(i) < b \).
  Equivalently,
  \begin{equation} \label{eq:ti}
    t_i = |\{j \in [n]: j < \sigma(i), \sigma(j) \ge j\}|
    - |\{j \in [n]: j < \sigma(i), \sigma(j) < j\}| - 1.
  \end{equation}
  We define the weight of \( i \) by \( \wt_\sigma(i) = Aq^{t_i} \).
\item Finally, if a vertex \( i \) is a decorated ending point,
  let \( s_i \) be the number of starting points or fixed points
  that come after \( i \).
  That is,
  \begin{equation} \label{eq:si}
    s_i = |\{j \in [n]: j > i, \sigma(j) \ge j \}|.
  \end{equation}
  We define the weight of \( i \) by \( \wt_\sigma(i) = Cq^{2s_i} \).
\end{itemize}
For each decorated matching \( \sigma \),
we define its weight by
\[
  \wt(\sigma) = q^{\cro(\sigma)} \prod_{i = 1}^n \wt_\sigma(i),
\]
where \( \cro(\sigma) \) is the number of crossings of
the underlying undecorated matching.
We define \( \wt(\sigma) = 1 \) if \( \sigma \in M^*(0, 0) \).
\begin{exam}
  Consider the matching
  \begin{equation}\label{eq:example-matching}
    \sigma = (1,3)(2,5)(4) \in M(5,1).
  \end{equation}
  The two edges form a crossing, and the fixed point \( 4 \)
  forms a crossing with the edge \( (2,5) \).
  Thus \( \cro(\sigma) = 2 \).
  For the starting and ending points, we have
  \[
    t_1 = t_2 = 1, \qquad s_3 = 1, \qquad s_5 = 0.
  \]
  Indeed, the edge \( (2,5) \) contributes to \( t_1 \),
  while the fixed point \( 4 \) contributes to both \( t_2 \)
  and \( s_3 \).

  If precisely the vertices \( 1, 3, 5 \) are decorated,
  the resulting decorated matching has weight
  \[
    q^2(Aq)(Cq^2)(C) = AC^2q^5.
  \]
  \begin{samepage}
  This decoration is shown below, where the filled vertices are decorated.
  
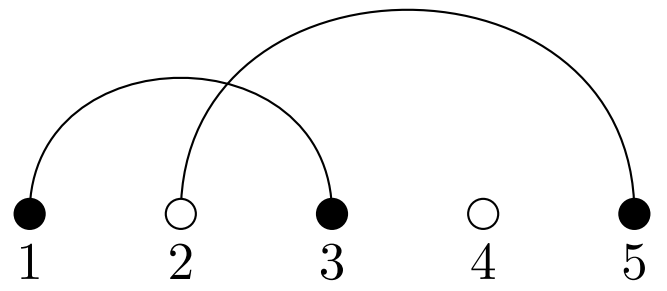
\begin{center}
    \begin{tikzpicture}[x=1cm,y=1cm]
      \draw (1,0) .. controls (1,1.2) and (3,1.2) .. (3,0);
      \draw (2,0) .. controls (2,1.8) and (5,1.8) .. (5,0);
      \foreach \i in {1,3,5}
        \node[circle,draw,fill=black,inner sep=2pt] at (\i,0) {};
      \foreach \i in {2,4}
        \node[circle,draw,fill=white,inner sep=2pt] at (\i,0) {};
      \foreach \i in {1,...,5}
        \node[below] at (\i,-0.08) {\( \i \)};
    \end{tikzpicture}
  \end{center}
  \end{samepage}

  Summing over all \( 2^4 = 16 \) decorations of this matching gives
  \[
    q^2(1+Aq)^2(1+Cq^2)(1+C).
  \]
\end{exam}

\begin{thm}\label{thm:comb}
  We have
  \[
    \tau_{n, k}(A, C; q) = \sum_{\sigma \in M^*(n, k)} \wt(\sigma).
  \]
\end{thm}
\begin{proof}
  We argue by induction on \( n \).
  The case \( n = 0 \) is immediate.

  Now suppose \( n \ge 1 \).
  For each \( \sigma \in M^{*}(n, k) \),
  there are two cases:
  \( n \) is a fixed point in the underlying matching,
  or \( n \) is not a fixed point in the underlying matching.

  \textbf{Case I: \( n \) is a fixed point.}
  Suppose that \( n \) is a fixed point in \( \sigma \in M^*(n, k) \).
  Then removing \( n \) gives a decorated matching
  \( \sigma' \in M^*(n-1, k-1) \).
  For each decorated ending point \( i \), removing the fixed point \( n \)
  changes its weight according to
  \[
    \wt_{\sigma'}(i) = q^{-2}\wt_{\sigma}(i).
  \]
  Therefore, we have
  \[
    \sum_{\substack{\sigma \in M^*(n, k) \\ \sigma(n) = n}} \wt(\sigma)
    = \sum_{\sigma' \in M^*(n-1, k-1)} \wt(\sigma')|_{C = Cq^2},
  \]
  which equals \( \tau_{n-1, k-1}(A, Cq^2; q) \) by the induction hypothesis.

  \textbf{Case II: \( n \) is not a fixed point.}
  Suppose that \( n \) is not a fixed point and \( \sigma(n) = j \)
  for \( \sigma \in M^*(n, k) \).
  Remove the edge \( (j, n) \), delete the vertex \( n \),
  and erase the decoration of \( j \), if any.
  This gives a decorated matching \( \sigma' \in M^*(n-1, k+1) \).

  Conversely, adding the vertex \( n \) and an edge \( (j, n) \)
  to \( \sigma' \in M^*(n-1, k+1) \) gives an element
  \( \sigma \in M^*(n, k) \).
  To add the edge \( (j, n) \), we must choose a
  fixed point \( j \) of \( \sigma' \).
  Therefore, there are \( k+1 \) choices of \( j \).
  If we choose the \( c \)th last fixed point,
  adding the edge \( (j, n) \) creates new crossings with fixed points
  between \( j \) and \( n \).
  Whenever an edge forms a crossing with the edge \( (j, n) \) in \( \sigma \),
  the edge forms a crossing with the fixed point \( j \) in \( \sigma' \).
  Therefore, the new edge \( (j, n) \) does not create additional crossings
  other than the crossings with fixed points.
  Thus the number of crossings increases by \( c-1 \),
  that is, \( \cro(\sigma) = \cro(\sigma') + c - 1 \).

  We may now decorate or leave undecorated each of the vertices
  \( j \) and \( n \).
  The weights \( \wt_\sigma(j) \) and \( \wt_\sigma(n) \)
  are determined according to the following choices.
  \begin{itemize}
  \item If \( j \) is left undecorated, then \( \wt_{\sigma}(j) = 1 =\wt_{\sigma'}(j)  \).
  \item If \( j \) is decorated, then
    \( \wt_{\sigma}(j) = Aq^{k} = Aq^{k}\wt_{\sigma'}(j) \).
  \item If \( n \) is left undecorated, then \( \wt_{\sigma}(n) = 1 \).
  \item If \( n \) is decorated, then \( \wt_{\sigma}(n) = C \).
  \end{itemize}
  Combining these choices, we obtain
  \[
    \sum_{\substack{\sigma \in M^*(n, k) \\ \sigma(n) < n}} \wt(\sigma)
    = \sum_{j = 1}^{n-1} \sum_{\substack{\sigma \in M^*(n, k) \\ \sigma(n) = j}} \wt(\sigma)
    = \sum_{\sigma' \in M^*(n-1, k+1)} \sum_{c=1}^{k+1} q^{c-1}
    (1 + Aq^{k})(1 + C)\wt(\sigma').
  \]
  By the induction hypothesis, this expression equals
  \( [k+1]_q(1 + Aq^k)(1 + C)\tau_{n-1, k+1}(A, C; q) \).

  By \Cref{thm:main}, we have
  \[
    \sum_{\sigma \in M^*(n, k)} \wt(\sigma)
    = \tau_{n-1, k-1}(A, Cq^2; q) + [k+1]_q(1 + Aq^k)(1 + C)\tau_{n-1, k+1}(A, C; q)
    = \tau_{n, k}(A, C; q).
  \]
\end{proof}

\begin{remark}
  Let \( 0\le k\le n \) with \( n-k \) even.
  Under the specialization \( A=C=0 \), we have
  \( \lambda^\nu_r(0,0;q)=[r]_q \), so
  \( p^\nu_r(x;0,0;q) \) is the monic continuous \( q \)-Hermite polynomial
  in the normalization of \cite[Eq.~(2.11)]{ISV}.

  Every matching with a decorated vertex has weight zero, while
  an undecorated matching has weight \( q^{\cro(\sigma)} \).
  Thus, \Cref{thm:comb} gives
  \[
    \tau_{n,k}(0,0;q)
    = \nu_{n,k}(0,0;q)
    = \sum_{\sigma\in M(n,k)} q^{\cro(\sigma)}.
  \]
  This is the partial matching formula for the mixed moments
  of continuous \( q \)-Hermite polynomials given by
  Cigler and Zeng~\cite{CiglerZeng2011}.
  Their exponent counts both crossings between edges
  and pairs consisting of an unmatched vertex and an edge covering it,
  so it equals \( \cro(\sigma) \) in our notation.
  In particular, \( k=0 \) recovers the classical matching formula
  for the even \( q \)-Hermite moments~\cite[Eq.~(3.8)]{ISV}.

  When \( A=C=q=1 \), every decorated matching has weight \( 1 \).
  Hence the sum of the coefficients of \( \tau_{n,k}(A,C;q) \) is
  \[
    \tau_{n,k}(1,1;1)
    = |M^*(n,k)|
    = 2^{n-k}\binom{n}{k}(n-k-1)(n-k-3)\cdots 1,
  \]
  where the product is \( 1 \) when \( n=k \).
\end{remark}

\begin{figure}
  \centering
  \newcommand{\matchingdiagram}[5]{%
    \begin{tikzpicture}[x=0.55cm,y=0.55cm,
      baseline=(current bounding box.center)]
      \draw (#1,0) .. controls (#1,1) and (#2,1) .. (#2,0);
      \foreach \v in {1,2,3}
        \node[circle,draw,fill=white,inner sep=2pt] at (\v,0) {};
      \ifnum#3=1
        \node[circle,draw,fill=black,inner sep=2pt] at (#1,0) {};
      \fi
      \ifnum#4=1
        \node[circle,draw,fill=black,inner sep=2pt] at (#2,0) {};
      \fi
      \foreach \v in {1,2,3}
        \node[below,font=\scriptsize] at (\v,-0.08) {\(\v\)};
      \node at (2,-1.1) {\(#5\)};
    \end{tikzpicture}%
  }
  \setlength{\tabcolsep}{10pt}
  \begin{tabular}{ccccc}
    \( (1,2)(3) \)
      & \matchingdiagram{1}{2}{0}{0}{1}
      & \matchingdiagram{1}{2}{1}{0}{A}
      & \matchingdiagram{1}{2}{0}{1}{Cq^2}
      & \matchingdiagram{1}{2}{1}{1}{ACq^2} \\[10pt]
    \( (1,3)(2) \)
      & \matchingdiagram{1}{3}{0}{0}{q}
      & \matchingdiagram{1}{3}{1}{0}{Aq^2}
      & \matchingdiagram{1}{3}{0}{1}{Cq}
      & \matchingdiagram{1}{3}{1}{1}{ACq^2} \\[10pt]
    \( (2,3)(1) \)
      & \matchingdiagram{2}{3}{0}{0}{1}
      & \matchingdiagram{2}{3}{1}{0}{Aq}
      & \matchingdiagram{2}{3}{0}{1}{C}
      & \matchingdiagram{2}{3}{1}{1}{ACq}
  \end{tabular}
  \caption{All twelve decorated matchings in \( M^*(3,1) \).
    The columns indicate which endpoints of the unique edge are decorated.
    Filled vertices are decorated, and the weight of each matching is shown below it.}
  \label{fig:decorated-matchings-three}
\end{figure}
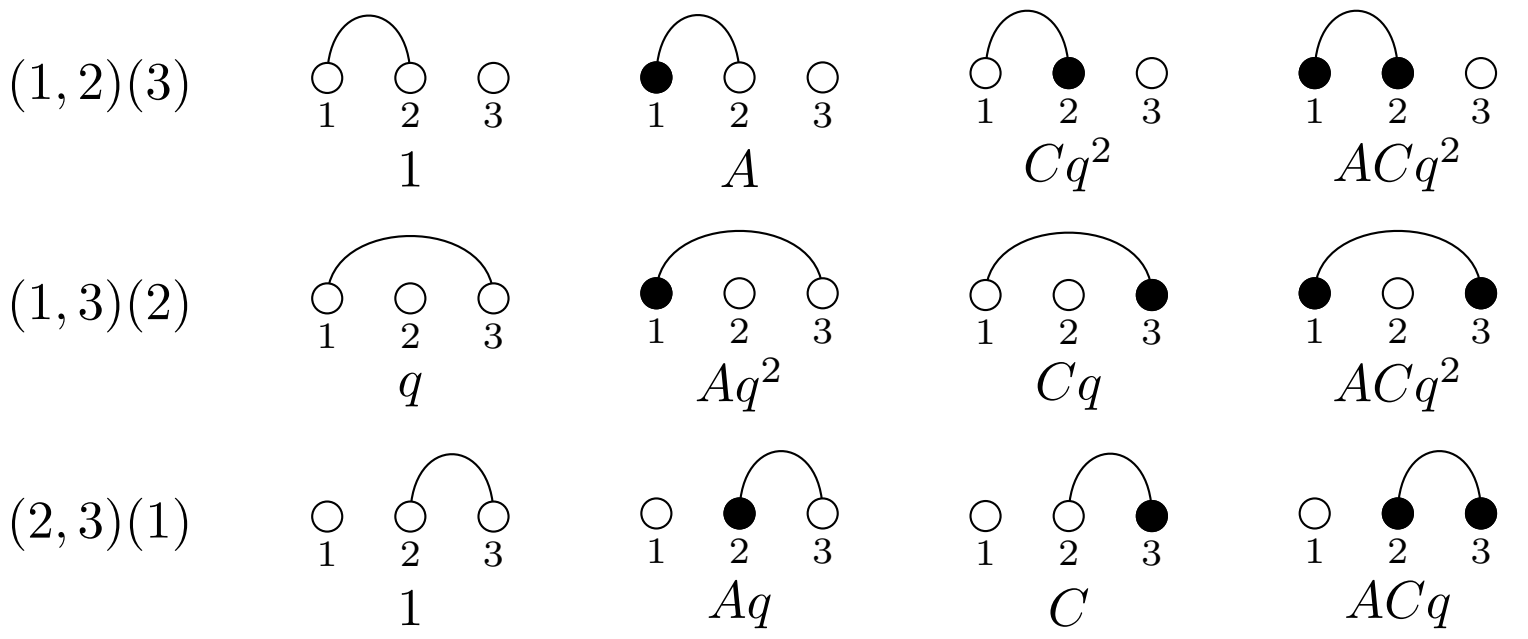

\begin{exam}
  The twelve decorated matchings in \( M^*(3,1) \) and their weights
  are shown in \Cref{fig:decorated-matchings-three}.
  There are three underlying matchings in \( M(3,1) \).
  Their contributions, after summing over all decorations, are
  \[
    \begin{array}{c|c|c|c|l}
      \sigma & \cro(\sigma) & t_i & s_j & \text{contribution} \\
      \hline
      (1,2)(3) & 0 & t_1=0 & s_2=1 & (1+A)(1+Cq^2) \\
      (1,3)(2) & 1 & t_1=1 & s_3=0 & q(1+Aq)(1+C) \\
      (2,3)(1) & 0 & t_2=1 & s_3=0 & (1+Aq)(1+C)
    \end{array}
  \]
  Here \( i \) and \( j \) are the starting point and ending point
  of the unique edge, respectively.
  Therefore, \Cref{thm:comb} gives
  \begin{align*}
    \tau_{3,1}(A,C; q)
    &= (1+A)(1+Cq^2) + (1+q)(1+Aq)(1+C) \\
    &= 2+q+(1+q+q^2)(A+C)+(q+2q^2)AC.
  \end{align*}
\end{exam}

Let \( 0 \le k \le n\), with \( n - k \) even,
and put \( m = \frac{n-k}{2} \).
It is known that the moments of the Askey--Wilson polynomials are symmetric
in \( a, b, c, d \)~\cite[Theorem~2.2]{AskeyWilson1985}.
Therefore, by definition, we have the symmetry
\[
  \tau_{n, k}(A, C; q) = \tau_{n, k}(C, A; q).
\]
We next turn to reciprocity. The recurrence coefficients satisfy
\( \lambda_n^\nu(A^{-1},C^{-1};q^{-1}) = -q \lambda_n^\nu(A,C;q) \).
Together with the definition of \( \tau_{n,k} \),
this gives the reciprocity relation
\[
  \tau_{n,k}(A,C;q)
  =A^mC^m q^{m(2k+m-1)}
  \tau_{n,k}(A^{-1},C^{-1};q^{-1}).
\]

The decorated matching model suggests the following problem.

\begin{problem}
  Find explicit bijections on \(M^*(n,k)\) that give
  direct combinatorial proofs of the symmetry and reciprocity relations above.
\end{problem}

\section{Multivariate generalization}
\label{sec:multi}
We now extend \Cref{thm:main} to multivariate moments.
Recall the definition of \( \{p^\nu_n(x)\}_{n \ge 0} \) in \eqref{eq:recurrence-nu}.
Consider its multivariate lift \( \{p^\nu_{\lambda}(\x_n)\}_{\lambda \in \Par_n} \).

Let 
\[
  V_{\lambda, \rho} = V_{\lambda, \rho}(A, C; q) :=
  \det (\nu_{\lambda_i + n - i, \rho_j + n - j}(A, C; q))_{i, j = 1}^n
\]
be the multivariate mixed moment of \( p^\nu_{\lambda}(\x_n) \),
and define \( T_{\lambda, \rho} = T_{\lambda, \rho}(A, C; q) \) to be
the minimal numerator of \( V_{\lambda, \rho} \).
If \( \rho = \emptyset \), we write
\[
  V_{\lambda} = V_{\lambda, \emptyset}, \qquad T_{\lambda} = T_{\lambda, \emptyset}.
\]
The following theorem gives a multivariate extension of the Kim--Stanton conjecture.

\begin{thm}\label{thm:multi}
  Let \( \lambda = (\lambda_1 ,\dots, \lambda_n) \in \Par_n \). Then
  \[
    (-1)^{\lambda_2 + \lambda_4 + \cdots } T_{\lambda}(A, C; q)
  \]
  is a polynomial in \( A, C, q \) with nonnegative integer coefficients.
\end{thm}

\begin{exam}
  Let \( n=3 \) and \( \lambda=(2,1,1) \).
  Suppressing the parameters \( (A,C;q) \), we obtain
  \[
    V_{(2,1,1)}=
    \begin{vmatrix}
      \nu_{4,2} & 0 & \nu_{4,0}\\
      1 & 0 & \nu_{2,0}\\
      0 & 1 & 0
    \end{vmatrix}
    =
    -\frac{[3]_q(1+A)(1+C)(1+Aq^2)(1+Cq^2)}
    {(1-ACq^3)(1-ACq^5)}.
  \]
  Expanding the determinant gives \( \nu_{4,0}-\nu_{4,2}\nu_{2,0} \).
  The factor \( 1-ACq \), which occurs in the denominators of both
  terms, cancels in their difference.
  The numerator and the denominator are coprime, and the denominator
  agrees with the product in \Cref{prop:deno}.
  Since \( (-1)^{\lambda_2}=-1 \), we obtain
  \[
    -T_{(2,1,1)}(A,C;q)
    =[3]_q(1+A)(1+C)(1+Aq^2)(1+Cq^2),
  \]
  which has nonnegative integer coefficients, as asserted by
  \Cref{thm:multi}.
\end{exam}

\Cref{thm:multi} does not extend to arbitrary multivariate mixed moments.
For \( \lambda = (7, 4) \) and \( \rho = (3, 0) \),
\( T_{\lambda, \rho}(A, C; q) \)
has constant term \( 26 \) and coefficient \( -13 \) on the monomial \( AC^2q^{19} \).
Thus the minimal numerator \( T_{\lambda, \rho}(A, C; q) \) does not have nonnegative coefficients
up to sign.

\begin{remark}
  The specialization \( q = 1 \) produces additional cancellations
  that do not occur in the general case.
  The minimal numerator of \( V_{\lambda, \rho}(A, C; 1) \)
  has nonnegative coefficients up to an overall sign.
\end{remark}

We now determine the minimal numerator explicitly.
If \( |\lambda| \) is odd, then \( V_\lambda=0 \),
so its minimal denominator and numerator are \( 1 \) and \( 0 \),
respectively.
For the denominator argument below, we may therefore assume
that \( |\lambda| \) is even.
We introduce a denominator of \( V_\lambda \) and show that
it is minimal whenever \( V_\lambda\neq0 \).

\begin{prop}\label{prop:deno}
  For \( \lambda \in \Par_n \), let
  \[
    D_\lambda = \prod_{i=1}^n \prod_{\substack{1 \le j \le \lambda_i \\ i - j \text{ odd}}}
    (1 - ACq^{2n-2-i+j}).
  \]
  Then \( D_\lambda V_\lambda(A, C; q) \) is a polynomial in \( A, C, q \).
\end{prop}
\begin{proof}
  First, the minimal denominators of \( \nu_{n, k} \) are of the form
  \( \prod_{i \in I}(1 - ACq^{i}) \) for some multiset \( I \)
  of odd integers.
  Therefore, the minimal denominator of \( V_\lambda \) has the same form.
  
  By the definition of \( V_\lambda \), we have
  \begin{align*}
    V_{\lambda}
    &= \det(\nu_{\lambda_i + n - i, n - j})_{i,j=1}^n \\
    &= \det\left(\left(\frac{2}{\sqrt{1-q}}\right)^{\lambda_i - i + j} \times
      \mu_{\lambda_i + n - i, n - j}(a, -a, c, -c; q)\right)_{i,j=1}^n \\
    &= \left( \frac{2}{\sqrt{1-q}} \right)^{|\lambda|} M_{\lambda}^{AW}(a, -a, c, -c; q),
  \end{align*}
  where \( M_{\lambda}^{AW} \) is the multivariate moment of the Askey--Wilson polynomials.
  By \cite[Theorem~4.12]{Cho2026a}, 
  \begin{align*}
    M^{AW}_\lambda(a, -a, c, -c; q)
    &= \sum_{\rho \subseteq \lambda} (-2a)^{-|\rho|} H_{\lambda , \rho}(f)
      s_\rho(1, q, \dots, q^{n-1}) \\
    &\quad \times \prod_{i=1}^n q^{\binom{n-i}{2}-\binom{\rho_i+n-i}{2}} \frac
      {(-A q^{n-i}; q)_{\rho_i} (AC q^{2n-2i};q^2)_{\rho_i}}
      {\left(AC q^{2 n-i-1};q \right)_{\rho_i}},
  \end{align*}
  where \( f = (f_0,f_1,\dots) \) is the sequence given by
  \( f_k = (aq^k + a^{-1}q^{-k})/2 \), and
  \[
    H_{\lambda, \rho}(f) :=
    \det(h_{\lambda_i - i - \rho_j + j}(f_0 ,\dots, f_{\rho_j + n - j}))_{i,j=1}^n.
  \]
  Here \( h_m(z_0 ,\dots, z_s) \) is the complete homogeneous symmetric polynomial,
  with
  \[
    h_0(z_0 ,\dots, z_s) = 1, \qquad h_m(z_0 ,\dots, z_s) = 0 \text{ for } m < 0.
  \]
  Therefore, the minimal denominator of \( V_\lambda \) divides
  \begin{align*}
    &(A(1-q))^{|\lambda|/2}q^{\sum_{i=1}^n\left(\binom{\lambda_i+n-i}{2} - \binom{n-i}{2}\right)}
    \prod_{i=1}^n\left(AC q^{2 n-i-1};q \right)_{\lambda_i} \\
    &= (A(1-q))^{|\lambda|/2}q^{\sum_{i=1}^n\left(\binom{\lambda_i+n-i}{2} - \binom{n-i}{2}\right)}
      \prod_{i=1}^n\prod_{1 \le j \le \lambda_i}(1 - AC q^{2n-2-i+j}).
  \end{align*}

  As shown above, the minimal denominator is
  of the form \( \prod_{i \in I}(1 - ACq^{i}) \)
  for some multiset \( I \) consisting of odd integers.
  Therefore, the minimal denominator divides
  \[
    D_\lambda = \prod_{i=1}^n \prod_{\substack{1 \le j \le \lambda_i \\ i - j \text{ odd}}}
    (1 - ACq^{2n-2-i+j}).
  \]
  Thus, \( D_\lambda V_\lambda \) is a polynomial in \( A, C, q \).
\end{proof}

To prove \Cref{thm:multi}, we first establish the following lemmas.
Define
\[
  \delta_j = \begin{cases}
    \frac{j - 2}{2} & \mbox{if \( j \) is even}, \\
    \frac{j - 1}{2} & \mbox{if \( j \) is odd}. \\
  \end{cases}
\]
\begin{lem}\label{lem:shift}
  Let
  \[
    V'_{\lambda}(A, C; q) = \det(\nu_{\lambda_i+n-i, n-j}(A, Cq^{2\delta_j}; q))_{i, j=1}^n.
  \]
  Then 
  \[
    V_\lambda(A, C; q) = V'_\lambda(A, C; q).
  \]
\end{lem}
\begin{proof}
  Let \( M_\lambda \) and \( M'_\lambda \)
  be the matrices given by
  \begin{align*}
    M_\lambda &= (\nu_{\lambda_i+n-i, n-j}(A, C; q))_{i, j=1}^n, \\
    M'_\lambda &= (\nu_{\lambda_i+n-i, n-j}(A, Cq^{2\delta_j}; q))_{i, j=1}^n,
  \end{align*}
  so that \( \det M_\lambda = V_\lambda \) and \( \det M'_\lambda = V'_\lambda \).
  
  Equation \eqref{eq:1} gives
  \[
    \nu_{n, k}(A, Cq^2; q) = \nu_{n, k}(A, C; q) - \eta_k \nu_{n, k+2}(A, C; q),
  \]
  where 
  \[
    \eta_k = \frac{C(1+Aq^k)(1+Aq^{k+1})(1-q)[k+1]_q [k+2]_q}{(1-ACq^{2k+1})(1-ACq^{2k+3})}.
  \]
  Applying \eqref{eq:1} \( m \) times to \( \nu_{n,k}(A, Cq^{2m}; q) \) expresses
  \( \nu_{n,k}(A, Cq^{2m}; q) \) as a linear combination of
  \( \nu_{n, k + 2s}(A, C; q) \) for \( 0 \le s \le m \),
  with coefficients independent of \( n \).

  Therefore,
  \begin{align*}
    M'_\lambda
    &= (\nu_{\lambda_i+n-i, n-j}(A, Cq^{2\delta_j}; q))_{i, j=1}^n \\
    &= \left(
      \sum_{s = 0}^{\delta_j} c_{j, s}
      \nu_{\lambda_i+n-i, n-j+2s}(A, C; q) \right)_{i, j=1}^n,
  \end{align*}
  for coefficients \( c_{j, s} \) independent of \( i \).
  The term indexed by \( n-j+2s \) corresponds to column \( j-2s \) of \( M_\lambda \).
  Each shift preserves the coefficient of \( \nu_{n,k} \), so \( c_{j,0}=1 \).
  Thus, the \( j \)th column of \( M'_\lambda \) is a linear combination of columns
  \( j, j-2 ,\dots,j-2\delta_j \) of \( M_\lambda \).
  Therefore,
  \[
    \det(M'_\lambda) = \det(M_\lambda).
  \]
\end{proof}

Define
\[
  N_\lambda = (\tau_{\lambda_i+n-i, n-j}(A, Cq^{2\delta_j}; q))_{i,j=1}^n
\]
and 
\[
  T'_\lambda(A, C; q) = \det(N_\lambda).
\]

Regard the rows and columns of \( N_\lambda \) as labelled by
\( r_i = \lambda_i + n - i \) and \( n - j \), respectively.
Since \( \tau_{n, k} = 0 \) unless \( n - k \) is even,
an entry can be nonzero only if the row and column labels have the same parity.
Therefore, if \( T'_\lambda(A, C; q) \neq 0 \),
there is a unique matrix \( \widetilde{N}_\lambda \)
obtained by rearranging rows of \( N_\lambda \)
so that the \( (i, j) \)-th entry of \( \widetilde{N}_\lambda \) is zero unless \( i - j \) is even,
and the row labels decrease within each parity class.

\begin{lem}
  Assume that \( T'_\lambda(A, C; q) \neq 0 \) and
  let \( \varepsilon_\lambda = (-1)^{\lambda_2 + \lambda_4 + \cdots} \).
  Then we have
  \[
    \varepsilon_\lambda \det \widetilde{N}_\lambda = T'_\lambda(A, C; q).
  \]
\end{lem}
\begin{proof}
  Let \( r_i = \lambda_i + n - i \) as before.
  We count how many row exchanges we need to obtain \( \widetilde{N}_\lambda \)
  from \( N_\lambda \).
  Note that \( \tau_{n, k} \) is nonzero only if \( n - k \) is even.
  Therefore, the \( (i, j) \)-th entry of \( N_\lambda \)
  is nonzero only if \( r_i - n + j \) is even.
  In particular, since \( r_i = \lambda_i + n - i \),
  \( r_i \) and \( n-i \) have the same parity
  if and only if \( \lambda_i \) is even.
  Since the \( (i, i) \)-th entry of \( \widetilde{N}_\lambda \) is nonzero,
  the rows in \( N_\lambda \) corresponding to odd \( \lambda_i \) must be moved.
  In order to resolve the parity mismatch of \( r_i \) and \( n-i \),
  we need to exchange rows in even positions with rows in odd positions.
  If \( T'_\lambda \neq 0 \), the number of odd \( \lambda_i \) with odd \( i \)
  equals the number of odd \( \lambda_i \) with even \( i \).
  Pair the rows in these two groups and exchange the rows in each pair.
  Then, exchanging rows gives the sign \( (-1)^{\lambda_2 + \lambda_4 + \cdots} \).
  
  We next show that the sign does not change after
  rearranging \( r_i \) in decreasing order
  within the same parity classes.
  Therefore, we need to compute the parity of the number of inversions
  within the same parity classes.
  To compute this, we first see the number of inversions in the original row order,
  and track how the exchanges in the previous step change the inversions.
  Initially, the row labels are in decreasing order.
  Therefore, every pair of row labels of the same parity is an inversion,
  that is, the larger label comes before the smaller label.
  If \( n \) is even, we exchange
  a row in an even position having an odd label
  with a row in an odd position with an even label.
  If \( n \) is odd, we exchange
  a row in an even position having an even label
  with a row in an odd position with an odd label.
  In either case, both the labels and the positions of
  the two exchanged rows have opposite parity.
  Since the labels of the two rows have different parity,
  the label of each row between them has the same parity
  as exactly one of the two exchanged labels.
  The exchange reverses the inversion status of this pair.
  Furthermore, there are an even number of rows between them.
  The total number of inversions within the same parity classes
  changes by an even number.
  The rearrangement of \( r_i \) in decreasing order makes
  every pair of the same parity be an inversion,
  and therefore needs an even number of transpositions.
  Hence, the sign does not change after
  rearranging \( r_i \) in decreasing order
  within the same parity classes.
\end{proof}

\begin{exam}
  Let \( n=5 \) and \( \lambda=(3,3,2,1,1) \).
  Suppressing the parameters \( (A,Cq^{2\delta_j};q) \) of the entries in column \( j \), we have
  \begingroup
  \def\quad{\kern3.5pt}
  \newcommand{\matrixlabel}[1]{\text{\scriptsize\sffamily\bfseries #1}}
  \[
  \bordermatrix{
     & \matrixlabel{4} & \matrixlabel{3} & \matrixlabel{2} & \matrixlabel{1} & \matrixlabel{0} \cr
    \matrixlabel{7} & 0 & \tau_{7,3} & 0 & \tau_{7,1} & 0 \cr
    \matrixlabel{6} & \tau_{6,4} & 0 & \tau_{6,2} & 0 & \tau_{6,0} \cr
    \matrixlabel{4} & \tau_{4,4} & 0 & \tau_{4,2} & 0 & \tau_{4,0} \cr
    \matrixlabel{2} & 0 & 0 & \tau_{2,2} & 0 & \tau_{2,0} \cr
    \matrixlabel{1} & 0 & 0 & 0 & \tau_{1,1} & 0 \cr
  }
  \xrightarrow{\substack{(1,4)\\(2,5)}}
  \bordermatrix{
     & \matrixlabel{4} & \matrixlabel{3} & \matrixlabel{2} & \matrixlabel{1} & \matrixlabel{0} \cr
    \matrixlabel{2} & 0 & 0 & \tau_{2,2} & 0 & \tau_{2,0} \cr
    \matrixlabel{1} & 0 & 0 & 0 & \tau_{1,1} & 0 \cr
    \matrixlabel{4} & \tau_{4,4} & 0 & \tau_{4,2} & 0 & \tau_{4,0} \cr
    \matrixlabel{7} & 0 & \tau_{7,3} & 0 & \tau_{7,1} & 0 \cr
    \matrixlabel{6} & \tau_{6,4} & 0 & \tau_{6,2} & 0 & \tau_{6,0} \cr
  }
  \xrightarrow{\substack{(1,5)\\(2,4)}}
  \bordermatrix{
     & \matrixlabel{4} & \matrixlabel{3} & \matrixlabel{2} & \matrixlabel{1} & \matrixlabel{0} \cr
    \matrixlabel{6} & \tau_{6,4} & 0 & \tau_{6,2} & 0 & \tau_{6,0} \cr
    \matrixlabel{7} & 0 & \tau_{7,3} & 0 & \tau_{7,1} & 0 \cr
    \matrixlabel{4} & \tau_{4,4} & 0 & \tau_{4,2} & 0 & \tau_{4,0} \cr
    \matrixlabel{1} & 0 & 0 & 0 & \tau_{1,1} & 0 \cr
    \matrixlabel{2} & 0 & 0 & \tau_{2,2} & 0 & \tau_{2,0} \cr
  }.
  \]
  \endgroup
  The small bold numbers to the left of and above each matrix
  are its row and column labels, respectively.
  The first matrix is \( N_\lambda \) and the last matrix is \( \widetilde{N}_\lambda \).
  Each pair above an arrow gives the positions of the rows exchanged,
  in the indicated order.
  The first arrow matches the required parities,
  while the second sorts the labels within each parity class.
  The total numbers of inversions of labels within the same parity classes
  are \( 4,0,4 \), respectively,
  and each arrow uses two exchanges.
  Thus the total sign is \( +1=(-1)^{\lambda_2+\lambda_4} \).
\end{exam}

Recall from \Cref{cor:min-deno} that \( \tau_{n,k}(A,C;q) \)
is the minimal numerator of \( \nu_{n,k}(A,C;q) \).
However, a direct determinantal form
\[
  \det(\tau_{\lambda_i + n - i, n - j}(A, C; q))_{i, j = 1}^n
\]
need not be a numerator of  \( V_\lambda \).

The column shifts \( C\mapsto Cq^{2\delta_j} \) in the definition
of \( T'_\lambda \) resolve this difficulty.
By \Cref{lem:shift}, these shifts leave \( V_\lambda \) unchanged,
and the following lemma shows that \( T'_\lambda \)
is a numerator of \( V_\lambda \).

\begin{lem}\label{lem:Tprime} 
  We have
  \[
    T'_\lambda(A, C; q) = D_\lambda V_\lambda(A, C; q),
  \]
  where \( D_\lambda \) is the polynomial defined in \Cref{prop:deno}.
\end{lem}
\begin{proof}
  For \( 0 \le k \le n \) with \( n-k \) even, the definition of
  \( \tau_{n, k}(A, C; q) \) gives
  \[
    \tau_{n, k}(A, C; q) = (ACq^{2k+1}; q^2)_{\frac{n-k}{2}}\nu_{n, k}(A, C; q).
  \]
  Put
  \[
    \alpha_i = \left\lfloor\frac{\lambda_i+2n-i-1}{2}\right\rfloor,
    \qquad
    \beta_j = n-1-\left\lfloor\frac{j}{2}\right\rfloor.
  \]
  Whenever \( \lambda_i-i+j \) is a nonnegative even integer, we have
  \[
    (ACq^{2\delta_j+2n-2j+1}; q^2)_{\frac{\lambda_i-i+j}{2}}
    = \frac{(ACq; q^2)_{\alpha_i}}{(ACq; q^2)_{\beta_j}}.
  \]
  Therefore, we obtain
  \begin{align*}
    T'_\lambda(A, C; q)
    &= \det(\tau_{\lambda_i+n-i, n-j}(A, Cq^{2\delta_j}; q))_{i,j=1}^n \\
    &= \det\left(\frac{(ACq; q^2)_{\alpha_i}}{(ACq; q^2)_{\beta_j}}
      \nu_{\lambda_i+n-i, n-j}(A, Cq^{2\delta_j}; q)\right)_{i,j=1}^n.
  \end{align*}
  Factoring \( (ACq; q^2)_{\alpha_i} \) out of row \( i \) and
  \( (ACq; q^2)_{\beta_j}^{-1} \) out of column \( j \), we obtain
  \[
    T'_\lambda(A, C; q) = V'_\lambda(A, C; q)
    \frac{\prod_{i=1}^n(ACq; q^2)_{\alpha_i}}
    {\prod_{j=1}^n(ACq; q^2)_{\beta_j}}.
  \]

  We have
  \[
    \frac{\prod_{i=1}^n(ACq; q^2)_{\alpha_i}}{\prod_{j=1}^n(ACq; q^2)_{\beta_j}}
    = \prod_{i=1}^n \frac{(ACq; q^2)_{\alpha_i}}{(ACq; q^2)_{\beta_i}}
    = \prod_{i=1}^n \prod_{\beta_i < t \le \alpha_i} (1 - ACq^{2t - 1})
    = \prod_{i=1}^n \prod_{\substack{1 \le j \le \lambda_i \\ i-j \text{ is odd}}} (1 - ACq^{2n - 2 - i + j}).
  \]
  Therefore, we have
  \[
    T'_\lambda(A, C; q) = D_\lambda V_\lambda(A, C; q).
  \]
\end{proof}

Now we are ready to prove \Cref{thm:multi}.

\begin{figure}
  \centering
  \begin{tikzpicture}[x=0.66cm,y=0.66cm,font=\small,
    line cap=round,line join=round]
    \begin{scope}
      \node at (-3,4.05) {Before the swap};
      \draw[gray!35] (-6.2,0) -- (0.3,0);
      \draw[thick,orange!85!black,dashed]
        (-6,0) -- (-5,1) -- (-4,2) -- (-3,3)
        -- (-2,2) -- (-1,1) -- (0,0);
      \draw[thick,blue!70!black]
        (-4,0) -- (-3,1) -- (-2,0) -- (-1,1) -- (0,2);
      \foreach \point in {(-6,0),(-4,0),(-1,1),(0,2),(0,0)}
        \fill \point circle (1.5pt);
      \node[below] at (-6,0) {\(S_{\sigma(i+1)}\)};
      \node[below] at (-4,0) {\(S_{\sigma(i)}\)};
      \node[right] at (0,2) {\(E_i\)};
      \node[right] at (0,0) {\(E_{i+1}\)};
      \node at (-1.05,0.5) {\(v\)};
      \node[orange!85!black] at (-3,3.45) {\(p_{i+1}\)};
      \node[blue!70!black] at (-3,1.5) {\(p_i\)};
    \end{scope}
    \node at (1.8,1.5) {\(\longleftrightarrow\)};
    \begin{scope}[xshift=6.4cm]
      \node at (-3,4.05) {After the swap};
      \draw[gray!35] (-6.2,0) -- (0.3,0);
      \draw[thick,orange!85!black,dashed]
        (-6,0) -- (-5,1) -- (-4,2) -- (-3,3)
        -- (-2,2) -- (-1,1) -- (0,2);
      \draw[thick,blue!70!black]
        (-4,0) -- (-3,1) -- (-2,0) -- (-1,1) -- (0,0);
      \foreach \point in {(-6,0),(-4,0),(-1,1),(0,2),(0,0)}
        \fill \point circle (1.5pt);
      \node[below] at (-6,0) {\(S_{\sigma(i+1)}\)};
      \node[below] at (-4,0) {\(S_{\sigma(i)}\)};
      \node[right] at (0,2) {\(E_i\)};
      \node[right] at (0,0) {\(E_{i+1}\)};
      \node at (-1.05,0.5) {\(v\)};
      \node[orange!85!black] at (-3,3.45) {\(p'_i\)};
      \node[blue!70!black] at (-3,1.5) {\(p'_{i+1}\)};
    \end{scope}
  \end{tikzpicture}
  \caption{Tail swapping at the selected intersection \(v=(-k,l)\).
    The paths are relabelled by their ending points,
    so \(p'_i=p_{i+1,1}p_{i,2}\) and \(p'_{i+1}=p_{i,1}p_{i+1,2}\).}
  \label{fig:tail-switch}
\end{figure}
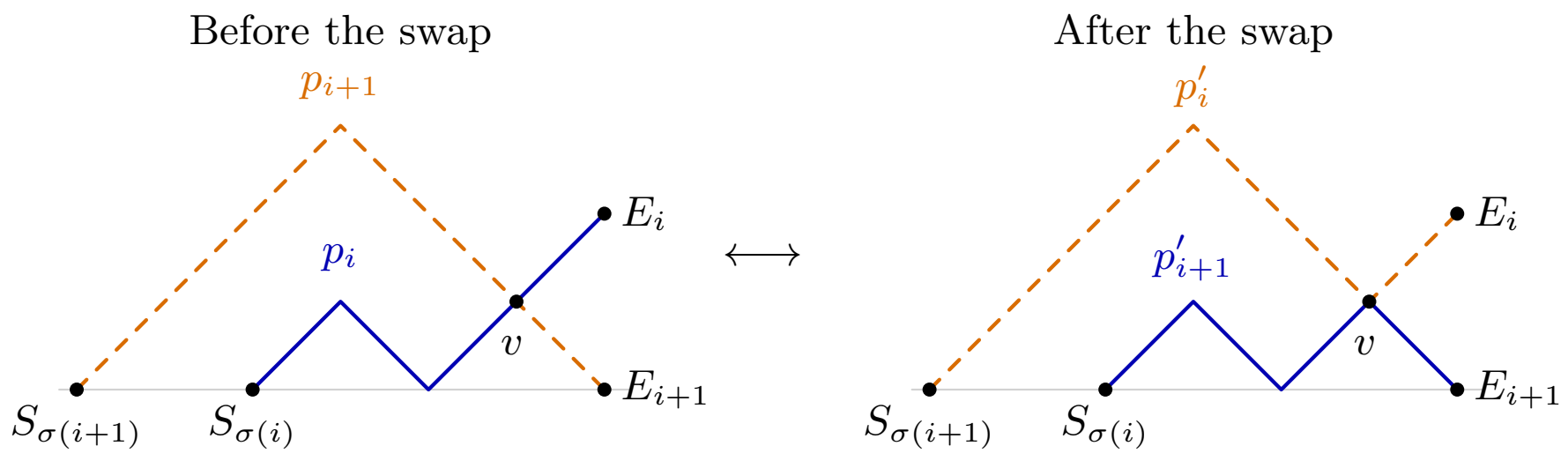

\begin{proof}[Proof of \Cref{thm:multi}]
  We first show the nonnegativity of \( \det \widetilde{N}_\lambda \)
  and prove the theorem by showing that \( T_\lambda(A, C; q) = T'_\lambda(A, C; q) \).
  The theorem is immediate if \( V_\lambda(A, C; q) = 0 \).
  Hence we assume that \( V_\lambda(A, C; q) \neq 0 \).
  
  \textbf{Step 1: Divide \( \widetilde{N}_\lambda \) into two blocks.}
  We assume that \( r_i \) are arranged so that
  \( r_i \) and \( n-i \) have the same parity.
  The nonzero entries in the even columns and odd columns
  form the blocks \( N_\lambda^e \) and \( N_\lambda^o \), respectively.
  \[
    N_\lambda^e
    =(\tau_{r_{2i},n-2j}(A,Cq^{2j-2};q))_{i,j=1}^{\lfloor n/2\rfloor}
  \]
  and
  \[
    N_\lambda^o
    =(\tau_{r_{2i-1},n-2j+1}(A,Cq^{2j-2};q))_{i,j=1}^{\lceil n/2\rceil}.
  \]
  Then \( \det(\widetilde{N}_\lambda) = \det(N_\lambda^e) \det(N_\lambda^o) \).

  It remains to show that both block determinants have nonnegative
  coefficients.
  We give the argument for \( \det(N_\lambda^e) \);
  the argument for \( \det(N_\lambda^o) \) is analogous.

  \textbf{Step 2: Construct a Dyck path model for multivariate moments.}
  Recall from \eqref{eq:path} that \( \tau_{n, k}(A, C; q) \) can be written as
  the weighted sum of Dyck paths from \( (-n, 0) \) to \( (0, k) \).
  Set \( m = \lfloor n/2 \rfloor \), and define the starting and ending points by
  \begin{align*}
    S_i &= (-r_{2i},0), & \mathcal S &= \{S_1,\dots,S_m\}, \\
    E_i &= (0,n-2i), & \mathcal E &= \{E_1,\dots,E_m\}.
  \end{align*}
  Each term in \( \det(N_\lambda^e) \) has the form
  \[
    \sgn(\sigma)\prod_{i=1}^m
    \tau_{r_{2\sigma(i)},n-2i}(A,Cq^{2i-2};q)
  \]
  for some \( \sigma\in S_m \).
  By \eqref{eq:path}, the \( i \)th factor \( \tau_{r_{2\sigma(i)},n-2i}(A,Cq^{2i-2};q) \)
  is the weighted sum of all Dyck paths
  \[
    p_i:S_{\sigma(i)}\longrightarrow E_i,
  \]
  with the weight of each path evaluated at \( C\mapsto Cq^{2i-2} \).

  Therefore, \( \det(N_\lambda^e) \) is the signed weighted sum over all tuples
  \( (p_1,\dots,p_m) \) of Dyck paths from \( \mathcal S \) to \( \mathcal E \).
  A tuple with \( p_i:S_{\sigma(i)}\to E_i \) contributes
  \[
    \sgn(\sigma)\prod_{i=1}^m
    \left.\wt(p_i)\right|_{C=Cq^{2i-2}}.
  \]
  In particular, the shift in the weight of \( p_i \) depends only on its
  ending point \( E_i \).

  \textbf{Step 3: Show that tail swapping is weight-preserving.}
  We prove that \( \det(N_\lambda^e) \) is the weighted sum over all
  nonintersecting tuples of Dyck paths from \( \{(-r_2, 0) ,\dots, (-r_{2m}, 0)\} \)
  to \( \{(0, n-2) ,\dots, (0, n-2m)\} \) by constructing
  a sign-reversing involution.

  Suppose that \( (p_1,\dots,p_m) \) is an intersecting tuple of Dyck paths from
  \[
    \{(-r_2,0),\dots,(-r_{2m},0)\}
    \quad\text{to}\quad
    \{(0,n-2),\dots,(0,n-2m)\}.
  \]
  Choose the rightmost intersection; if several occur at the same
  horizontal coordinate, choose the lowest one.
  This intersection must be shared by the adjacent paths
  \begin{align*}
    p_i &: (-r_{2\sigma(i)},0) \rightarrow (0,n-2i), \\
    p_{i+1} &: (-r_{2\sigma(i+1)},0) \rightarrow (0,n-2i-2),
  \end{align*}
  for some \( 1 \le i < m \) and \( \sigma \in S_m \).
  Denote the intersection point by \( (-k,l) \).
  Let the heads of these paths be
  \begin{align*}
    p_{i,1} &: (-r_{2\sigma(i)},0) \rightarrow (-k,l), \\
    p_{i+1,1} &: (-r_{2\sigma(i+1)},0) \rightarrow (-k,l),
  \end{align*}
  and let their tails be
  \begin{align*}
    p_{i,2} &: (-k,l) \rightarrow (0,n-2i), \\
    p_{i+1,2} &: (-k,l) \rightarrow (0,n-2i-2).
  \end{align*}
  Thus \( p_i = p_{i,1}p_{i,2} \) and \( p_{i+1} = p_{i+1,1}p_{i+1,2} \).

  The number of up steps in \( p_{i, 2} \) is
  \[
    u_1 := \frac{n + k - 2i - l}{2}
  \]
  and the number of up steps in \( p_{i+1, 2} \) is
  \[
    u_1 - 1 = \frac{n + k - 2i - l}{2} - 1.
  \]
  Let 
  \begin{align*}
    w_{i, 1}(A, C; q) &:= \wt(p_{i, 1}),& w_{i+1, 1}(A, C; q) &:= \wt(p_{i+1, 1}), \\
    w_{i, 2}(A, C; q) &:= \wt(p_{i, 2}),& w_{i+1, 2}(A, C; q) &:= \wt(p_{i+1, 2}). \\
  \end{align*}
  Then 
  \[
    \wt(p_i) = w_{i, 1}(A, Cq^{2u_1}; q)w_{i, 2}(A, C; q), \qquad
    \wt(p_{i+1}) = w_{i+1, 1}(A, Cq^{2u_1-2}; q)w_{i+1, 2}(A, C; q).
  \]
  Therefore, the contribution of the paths \( p_i \) and \( p_{i+1} \)
  to the term indexed by \( \sigma \) is
  \[
    w_{i, 1}(A, Cq^{2u_1+2i-2}; q)w_{i, 2}(A, Cq^{2i-2}; q)
    w_{i+1, 1}(A, Cq^{2u_1-2+2i}; q)w_{i+1, 2}(A, Cq^{2i}; q). 
  \]

  Now swap the tails, as shown in Figure~\ref{fig:tail-switch}, so that
  \[
    p'_i = p_{i+1, 1}p_{i, 2}, \qquad p'_{i+1} = p_{i, 1}p_{i+1, 2}.
  \]
  Then
  \[
    \wt(p'_i) = w_{i+1, 1}(A, Cq^{2u_1}; q)w_{i, 2}(A, C; q), \qquad
    \wt(p'_{i+1}) = w_{i, 1}(A, Cq^{2u_1-2}; q)w_{i+1, 2}(A, C; q),
  \]
  and the contribution of the paths \( p'_i \) and \( p'_{i+1} \)
  to the term indexed by \( \sigma'=\sigma \circ (i, i+1) \)
  in the expansion of \( \det(N_\lambda^e) \) is
  \[
    w_{i+1, 1}(A, Cq^{2u_1+2i-2}; q)w_{i, 2}(A, Cq^{2i-2}; q)
    w_{i, 1}(A, Cq^{2u_1-2+2i}; q)w_{i+1, 2}(A, Cq^{2i}; q). 
  \]
  
  Tail swapping preserves the rightmost intersection and
  we have \( \sgn(\sigma')=-\sgn(\sigma) \).
  Thus tail swapping defines a sign-reversing involution.
  After cancellation, only the nonintersecting tuples remain,
  proving the claim.

  Since the row labels decrease within each parity class,
  every nonintersecting tuple has the identity endpoint permutation:
  any other pairing would force two paths to share a vertex.
  Thus all surviving contributions have positive sign.
  The weights of Dyck paths have nonnegative coefficients, so
  \( \det(N_\lambda^e) \) has nonnegative coefficients.

  The same argument applies to \( \det(N_\lambda^o) \). Hence
  \[
    \varepsilon_\lambda T'_\lambda(A,C;q)
    = \det(\widetilde{N}_\lambda)
    =\det(N_\lambda^e)\det(N_\lambda^o)
  \]
  has nonnegative coefficients.

  \textbf{Step 4: Show the equality of \( T_\lambda \) and \( T_\lambda' \).}
  By \Cref{lem:Tprime}, we have
  \[
    T_\lambda'(A, C; q) = T_\lambda(A, C; q) \prod_{i \in I}(1 - ACq^i)
  \]
  for some multiset \( I \).
  Set \( A = C = q = 1 \).
  The left-hand side is nonzero unless \( T_\lambda'(A, C; q) = 0 \).
  Therefore, the multiset \( I \) must be empty,
  which implies
  \[
    T_\lambda'(A, C; q) = T_\lambda(A, C; q).
  \]
\end{proof}

We now record the minimality of the denominator \( D_\lambda \)
defined in \Cref{prop:deno}.

\begin{cor}
  Let \( \lambda \in \Par_n \) with \( V_{\lambda}(A, C; q) \neq 0 \).
  Then the minimal denominator of \( V_{\lambda}(A, C; q) \)
  is \( D_\lambda \).
\end{cor}
\begin{proof}
  The proof of \Cref{thm:multi} gives
  \[
    T_\lambda = T'_\lambda = D_\lambda V_\lambda.
  \]
  Since \( V_\lambda \neq 0 \) and \( T_\lambda \) is its
  minimal numerator, the claim follows.
\end{proof}

\begin{remark}
  After adjusting the indices, one can prove that for \( \lambda \in \Par_n \),
  \( T_{\lambda, (k^n)}(A, C; q) \) is a polynomial in \( A, C, q \)
  with nonnegative coefficients up to sign.
\end{remark}

The dual analogues of the multivariate mixed moments are
the multivariate coefficients:
the coefficients in the Schur expansions of the multivariate
orthogonal polynomials.
We conjecture that a similar phenomenon appears for the
multivariate coefficients.
The same argument cannot be applied to the multivariate coefficients,
because the head swapping is not weight-preserving.

\begin{conj}
  Let \( \lambda\in\Par_n \).
  Then the minimal numerator of the coefficient of \( s_\emptyset(\x_n)=1 \)
  in the Schur expansion of \( p^\nu_\lambda(\x_n) \)
  has nonnegative integer coefficients up to an overall sign.
\end{conj}

We have checked this conjecture for \( n \le 6 \) and
\( \lambda_1 \le 6 \).
  
\section*{Acknowledgement}
The author thanks Jang Soo Kim for introducing the problem
and for providing helpful advice on revising the paper.

\section*{AI declaration}
This paper was written with the assistance of AI tools (ChatGPT 6 Pro).
In particular, \Cref{thm:main}, \Cref{lem:1}, and \Cref{lem:shift}
were suggested by AI tools, and the author rewrote the corresponding proofs
based on these suggestions.
The rest of the mathematical content was developed by the author.
AI tools were also used to help draw diagrams, write Sage code,
identify relevant literature, and improve the exposition.

\bibliographystyle{abbrv}

\end{document}